\documentclass[12pt]{article}
\usepackage{microtype} \DisableLigatures{encoding = *, family = * }
\usepackage{subeqn,multirow}
\usepackage{graphicx}
\usepackage{amsfonts,amssymb}
\usepackage{natbib}
\bibpunct{(}{)}{;}{a}{,}{,}
\usepackage[bookmarksnumbered=true, pdfauthor={Wen-Long Jin}]{hyperref}

\newcommand{\commentout}[1]{}

\newcommand{\ba}{\begin{array}}
        \newcommand{\ea}{\end{array}}
\newcommand{\bc}{\begin{center}}
        \newcommand{\ec}{\end{center}}
\newcommand{\bdm}{\begin{displaymath}}
        \newcommand{\edm}{\end{displaymath}}
\newcommand{\bds} {\begin{description}}
        \newcommand{\eds} {\end{description}}
\newcommand{\ben}{\begin{enumerate}}
        \newcommand{\een}{\end{enumerate}}
\newcommand{\beq}{\begin{equation}}
        \newcommand{\eeq}{\end{equation}}
\newcommand{\bfg} {\begin{figure}}
        \newcommand{\efg} {\end{figure}}
\newcommand{\bi} {\begin {itemize}}
        \newcommand{\ei} {\end {itemize}}
\newcommand{\bpp}{\begin{pspicture}}
        \newcommand{\epp}{\end{pspicture}}
\newcommand{\bqn}{\begin{eqnarray}} 
        \newcommand{\eqn}{\end{eqnarray}}
\newcommand{\bqs}{\begin{eqnarray*}}
        \newcommand{\eqs}{\end{eqnarray*}}
\newcommand{\bsq}{\begin{subequations}}
        \newcommand{\esq}{\end{subequations}}
\newcommand{\bsl} {\begin{slide}[8.8in,6.7in]}
        \newcommand{\esl} {\end{slide}}
\newcommand{\bss} {\begin{slide*}[9.3in,6.7in]}
        \newcommand{\ess} {\end{slide*}}
\newcommand{\btb} {\begin {table}[h]}
        \newcommand{\etb} {\end {table}}
\newcommand{\bvb}{\begin{verbatim}}
        \newcommand{\evb}{\end{verbatim}}

\newcommand{\m}{\mbox}
\newcommand {\der}[2] {{\frac {\m {d} {#1}} {\m{d} {#2}}}}
\newcommand {\pd}[2] {{\frac {\partial {#1}} {\partial {#2}}}}

\newcommand{\cas}[1]{{{\left \{ \ba #1 \ea \right. }}}

\newcommand{\reff}[1] {{{Figure \ref {#1}}}}
\newcommand{\refe}[1] {{{(\ref {#1})}}}

\newtheorem{theorem}{Theorem}[section]
\newtheorem{lemma}[theorem]{Lemma}
\newtheorem{corollary}[theorem]{Corollary}

\def\pmb#1{\setbox0=\hbox{$#1$}%
   \kern-.025em\copy0\kern-\wd0
   \kern.05em\copy0\kern-\wd0
   \kern-.025em\raise.0433em\box0 }

\def\eop{{\hfill $\blacksquare$}}

\def\dt     {{\Delta t}}

\def\eop{{\hfill $\blacksquare$}}

\def\la {{{\lambda}}}

\def\L{{{\mathcal{L}}}}

\begin {document}
\title{Performance analysis and signal design for a stationary signalized ring road} 
\author{Wen-Long Jin \footnote{Department of Civil and Environmental Engineering, California Institute for Telecommunications and Information Technology, Institute of Transportation Studies, 4000 Anteater Instruction and Research Bldg, University of California, Irvine, CA, USA 92697-3600. Email: wjin@uci.edu. Corresponding author} and Yifeng Yu \footnote{Department of Mathematics, University of California, Irvine, CA, USA. Email: yifengy@uci.edu}}

\maketitle

\begin{abstract} Existing methods for traffic signal design are either too simplistic to capture realistic traffic characteristics or too complicated to be mathematically tractable. In this study, we attempts to fill the gap by presenting a new method based on the LWR model for performance analysis and signal design in a stationary signalized ring road. We first solve the link transmission model to obtain an equation for the boundary flow in stationary states, which are defined to be time-periodic solutions in both flow-rate and density with a period of the cycle length. We then derive an explicit macroscopic fundamental diagram (MFD), in which the average flow-rate in stationary states is a function of both traffic density and signal settings. Finally we present simple formulas for optimal cycle lengths under five levels of congestion with a start-up lost time. With numerical examples we verify our analytical results and discuss the existence of near-optimal cycle lengths. This study lays the foundation for future studies on performance analysis and signal design for more general urban networks based on the kinematic wave theory. 

\end{abstract}

\section{Introduction}

Traffic signals have been widely deployed to resolve conflicts among various traffic streams and improve safety of drivers and pedestrians at busy urban intersections. But signalized intersections are also major network bottlenecks, inducing stop-and-go traffic patterns, travel delays, and vehicle emissions. 
Many efforts have been devoted to mitigating the congestion effects of isolated and coordinated intersections by optimally designing phase sequences, cycle lengths, green splits, offsets, and other parameters of traffic signals  \citep{papageorgiou2005review}. 

Most of existing signal design methods employ two types of traffic flow models: aggregate delay and bandwidth formulas or traffic simulation models. 
In the first type of methods, the cycle length of a signal is selected to minimize vehicles' average delays according to Webster's formula, the allocation of the total effective green time in a cycle to different phases depends on their respective flow-rates, and then offsets at different intersections are determined by optimizing the bandwidth \citep{miller1963settings,gartner1975optimization,roess2010traffic}. Such an approach is straightforward for designing either pretimed or actuated signals but has serious limitations: first, Webster's formula is derived based on the assumption of random Poisson arrival patterns of vehicles, but in reality arrival patterns are regulated by other signals; second, Webster's formula only applies to under-saturated intersections without accounting for impacts of congested downstream links; third, delay formulas used at the design stage are usually different from those used at the analysis stage \citep{dion2004intersection}; finally, exact relationships between bandwidths and vehicles' delays are not clear. Therefore, methods based on Webster's and other delay formulas are too simplistic, since they cannot capture traffic waves on a link, interactions among intersections, queue spillback, or other realistic traffic phenomena in a signalized road network.
In the second type of methods, various traffic flow models are used to simulate realistic traffic dynamics, and optimal control problems are formulated to find best signal settings simultaneously subject to given demand patterns   \citep{gazis1963over,gazis1964optimum,dans1976optimal,improta1984control,papageorgiou1995integrated,park1999traffic,chang2000optimal,chang2004modeling}. However, such methods are too detailed to be amenable to mathematical analyses and computationally costly for studying large-scale networks. In summary, existing methods for traffic signal design are either too simplistic to capture realistic traffic characteristics or too complicated to be analytically solvable. We believe that this is major reason for the lack of ``a systematic theory (even) for a one-way arterial'' \citep{newell1989theory}.

The existence of a gap between methods based on delay formulas and those based on traffic simulation has motivated us to develop a new approach for performance analysis and signal design in signalized networks. Our daily experience suggests that traffic patterns in an urban road network are relatively stationary during peak periods; that is, the locations and durations of congestion are stable from day to day. 
The new approach builds on the assumption of the existence of such stationary states.  Furthermore, we attempt to derive the average flow-rate in stationary states as a function of signal settings within the framework of the LWR model \citep{lighthill1955lwr,richards1956lwr}. In such stationary networks, the average flow-rate can serve as a performance measure and the objective function to find optimal signal settings, since a larger flow-rate at the same density leads to lower average delays and more efficient operations. Therefore, such a method is both physically realistic and mathematically tractable.

\begin{figure} \bc
$\ba{c@{\hspace{1in}}c}
\includegraphics[height=2.5in]{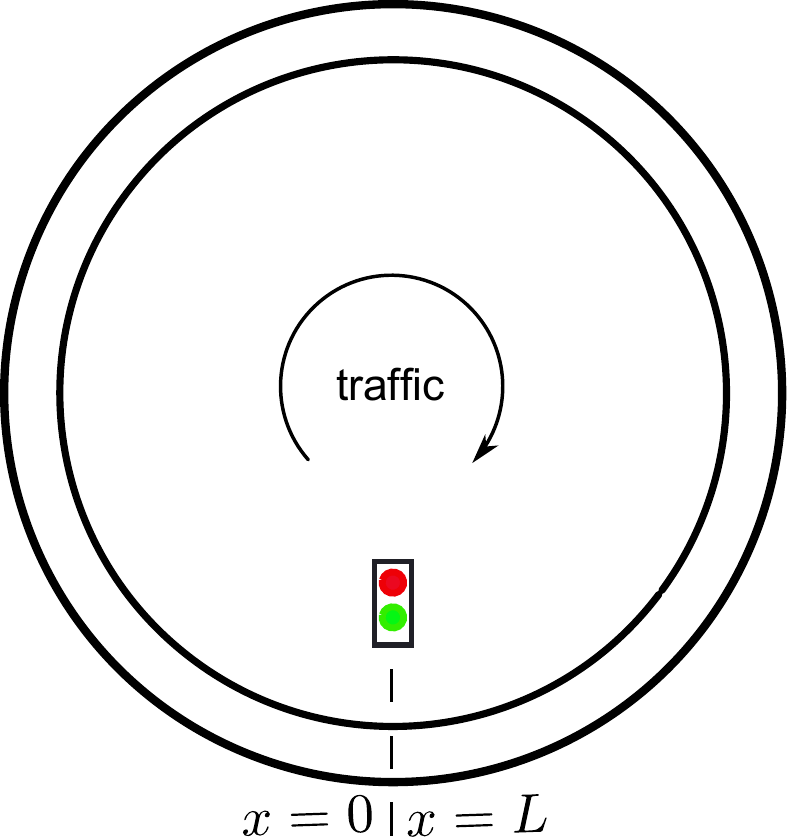} &
\includegraphics[height=2.5in]{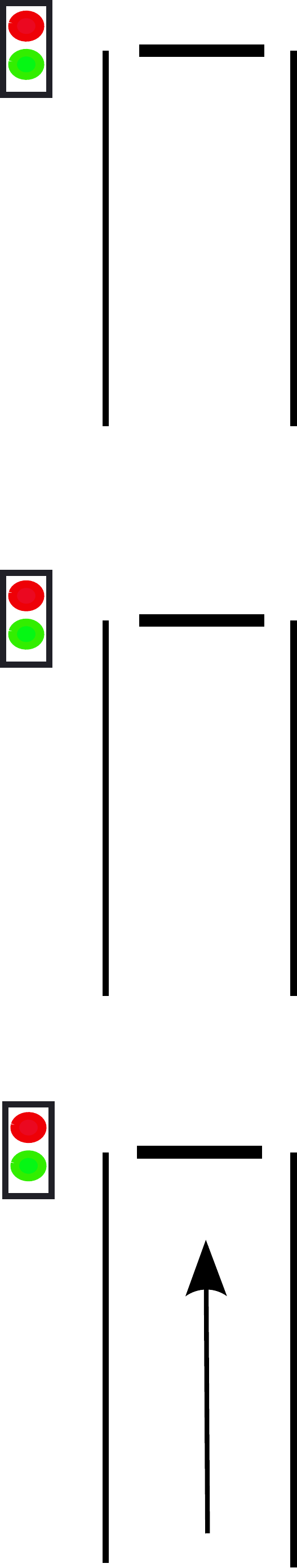} \\
{\mbox{\bf (a)}} &
    {\mbox{\bf (b)}}
\ea$
\caption{(a) A signalized ring road (b) An infinite street with identical roads} \label{signalized_ring_ill} \ec 
\end{figure}

In \citep{jin2014_hamiltonian}, it was proved that asymptotic periodic traffic patterns, which can be defined as stationary states, exist in a ring road with a pretimed signal, shown in \reff{signalized_ring_ill}(a). As a first step for developing a unified approach for signal design for general road networks, in this study we start with the signalized ring road, which is the simplest signalized network. The signalized ring road is equivalent to an infinite street without turning movements in \reff{signalized_ring_ill}(b), where all links, traffic conditions, and signal settings are identical. In this sense, the offset between two consecutive signals is 0, and signals follow a simultaneous progression in the northbound direction. In addition, we assume that each cycle has only two phases. Therefore, signal settings remain to be determined include the green split and the cycle length. 

In stationary urban road networks, it was postulated that there exits a relation between network-average flow and density in \citep{godfrey1969mechanism}. Such a relation is called the macroscopic fundamental diagram (MFD) and has been shown to be unique in homogeneous networks, but not in non-homogeneous ones with simulations and observations \citep{ardekani1987urban,mahmassani1987performance,olszewski1995area,geroliminis2008eus,buisson2009exploring,cassidy2011macroscopic,geroliminis2012effect}. In \citep{daganzo2007gridlock,geroliminis2013optimal}, regional demand control strategies were developed based on MFD. As a system-wide characteristic, MFD emerges from network traffic flow patterns, which are determined by network topology, signal and other control measures, and drivers' choices in destinations, modes, departure times, routes, and speeds. Some efforts have been devoted to deriving MFD in simple signalized networks from various traffic flow models. 
In \citep{gartner2004analysis}, with cellular automaton simulations for traffic on a ring road, which has  multiple identical signals, the relationship between flow-rate, density, and offset was obtained in relatively stationary states after a long time (2000 seconds), and it was found that offsets can have drastic impacts on the overall throughputs and, therefore, travel times on an arterial road.
 In \citep{daganzo2008analytical}, a variational method was proposed to calculate approximate MFD in a ring road with multiple signals, but no definitions of stationary states were provided, and impacts of signal settings were not considered. As far as we know, no simple guidelines for signal design were provided in this reference and its follow-up studies, even for a signalized ring road. 
In \citep{daganzo2011bifurcations}, MFD in a double-ring network with turning movements was studied with  heuristic double-bin approximations and cellular automaton simulations. 
In \citep{jin2013dring}, steady or stationary states in a signalized double-ring network were defined as asymptotically periodic traffic states within the framework of a network kinematic wave theory, and impacts of signal settings and turning movements on MFD in stationary states were simulated with Cell Transmission Model (CTM) \citep{daganzo1995ctm}. 
However, there has been no theory for the existence of such stationary states in general networks, and no explicit closed-form relation between signal settings and MFD is known, even for simple networks. In this study, we take one step further by deriving the average flow-rate in stationary states as a function of both density and signal settings, which can be used to find optimal signal settings at different congestion levels.

This study is enabled by the link transmission model (LTM) \citep{yperman2006mcl,yperman2007link}, which, together with Newell's simplified kinematic wave model \citep{newell1993sim}, is another formulation of the network kinematic wave theory based on the LWR model. In \citep{jin2015ltm}, two continuous formulations of the LTM were derived from the Hopf-Lax formula for the Hamilton-Jacobi equation of the LWR model. \footnote{Even though the LWR model was already applied to analyze the formation of queues at a signalized intersection in \citep{lighthill1955lwr}, no relationships between signal settings and system performance have been established.} For the signalized ring road, we first formulate and solve the LTM for boundary flows (Section 2), then derive an explicit formula for MFD in stationary states and analyze its relationship with signal cycle length (Section 3), and finally find optimal signal settings to maximize the average flow-rate in MFD (Section 4). In Section 5, we conclude the study with future directions.

\section{The link transmission model for a signalized ring road}
For a ring road with a length of $L$, as shown in \reff{signalized_ring_ill}(a), the $x$-axis increases in the traffic direction, and we place a signal at $x=0$ and $x=L$. We apply the LWR model to describe traffic dynamics in the signalized ring road. Due to the existence of the signal, the problem is much more challenging to solve within the traditional framework of hyperbolic conservation laws. Here we study it with the help of the LTM. 

\subsection{The LWR model and Hamilton-Jacobi equation}
The evolution of traffic density $k(x,t)$ on the ring road can be described by the LWR model: 
\bsq 
\bqn
\pd{k}t+\pd{B(x,t) Q(k)}x&=&0, \label{kw-signal}
\eqn
with a periodic boundary condition
\bqn
k(0,t)&=&k(L,t).
\eqn

Here we assume a triangular fundamental diagram \citep{munjal1971multilane,haberman1977model,newell1993sim}, 
\bqn
Q(k)&=&\min\{V k, (K-k)W\},
\eqn
where $V$ is the free-flow speed, $-W$ the shock wave speed in congested traffic, and $K$ the jam density. Thus the critical density is $\bar K=\frac{W}{V+W} K$, and the capacity $C=V \bar K$.

 The effect of the signal is captured by $B(x,t)=1-I(x) \cdot (1-\beta(t))$, where $I(x)$ determines the location of the traffic signal: $I(x)=\cas{{ll}1, &x=0, L\\0, &\m{otherwise}}$, $\pi\in (0,1)$ is the ratio of effective green time to the cycle length \footnote{Here we approximate a cycle, which comprises green, yellow, all red, and red intervals, by an effective green interval and an effective red interval.}, and the traffic light is effective green during $iT+[0,\pi T]$ and effective red during other time intervals:
\bqn
\beta(t)&=&\cas{{ll}1, &t-iT\in [0,\pi T], \quad i=0,1,2,\cdots\\0, & \m{otherwise}}
\eqn
\esq
where $T$ is the cycle length.

Following \citep{newell1993sim}, we can obtain the Hamilton-Jacobi equation of the LWR model, \refe{kw-signal}, by introducing a new state variable inside the spatial-temporal domain $\Omega=[0,L]\times[0,\infty)$, $A(x,t)$, which is the cumulative flow passing $x$ before $t$ and also known as a Moskowitz function \citep{moskowitz1965discussion}.
 Then $k=-A_x$, $q=A_t$, and the flow conservation equation is automatically satisfied since $A_{xt}=A_{tx}$. Furthermore the fundamental diagram and, therefore, the LWR model for the signalized ring road, \refe{kw-signal}, is equivalent to the following Hamilton-Jacobi equation \citep{evans1998pde}:
 \bsq \label{lwr_vt_signalized_ring}
\bqn
A_t-B(x,t)Q(-A_x)&=&0, \label{lwr_vt}
\eqn 
with a periodic boundary condition
\bqn
A_t(0,t)&=&A_t(L,t). \label{hj_periodic}
\eqn
\esq
Here the Hamiltonian is both space- and time-dependent, and $-Q(-A_x)$ is convex in $A_x$ for the triangular fundamental diagram. In \citep{jin2014_hamiltonian}, the solution of \refe{lwr_vt} was proved to be asymptotically periodic, and a relationship between the average density and flow-rate, i.e., MFD, was proved to exist, as the solution of a so-called ``cell problem''. However, the closed-form equation of MFD was not found.

\subsection{The link transmission model}
As shown in \reff{boundary_hopf_lax}, where the red bars denote the effective red light signals, we extend the spatial-temporal domain $\Omega$ periodically, since the signalized ring road is equivalent to a street with identical links.
We denote the initial condition by $N(x)=A(x,0)$ and  the average density on the ring road by $k_0$. Then $N(x)=N(x-L)-k_0L$. We further assume that the initial density is constant at $k_0$ on all links; i.e., initially $N(x)=-k_0 x$.
We denote the boundary flow by $G(t)=A(0,t)$ with $G(0)=0$ and the corresponding flow-rate by $g(t)$. Then
\bqn
\der{}{t} G(t)&=&g(t). \label{def:flow-rate}
\eqn
 Then from the periodic boundary condition \refe{hj_periodic}, we have $A(jL,t)=G(t)-jk_0L$ for $j=0,\pm1,\pm2,\cdots$.

In the following we present the continuous and discrete formulations of LTM, which can be used to solve $G(t)$ both analytically and numerically. Further from $G(t)$, we can solve \refe{lwr_vt} to obtain $A(x,t)$ for $x\in(0,L)$ by following Newell's minimization principle \citep{newell1993sim}, Daganzo's variational principle  \citep{daganzo2005variationalKW,daganzo2005variationalKW2}, the Hopf-Lax formula, optimal control principle, or the viscosity solution method \citep{evans1998pde}. But we are not concerned with traffic dynamics inside the ring road in this study.

\bfg\bc
\includegraphics[width=3in]{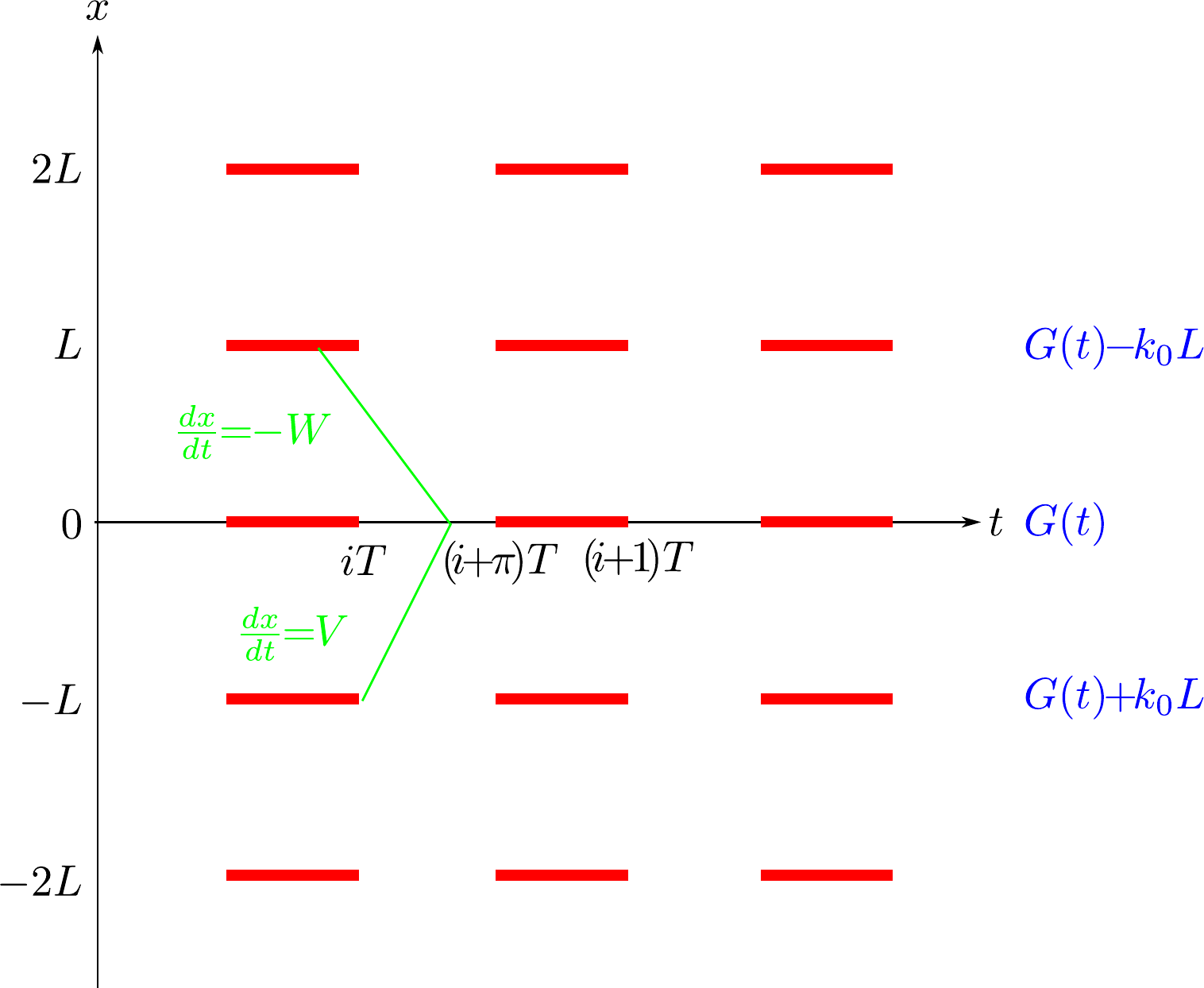} \caption{A periodic extension of the spatial-temporal domain $\Omega$} \label{boundary_hopf_lax}
\ec\efg

In the continuous formulation of LTM \citep{jin2015ltm}, we first define the demand, $d(t)$, of the link between $-L$ and 0:
\bsq\label{def:demand}
\bqn
d(t)&=&\cas{{ll} \min\left\{k_0 V+H(\la(t)), C\right\}, & t\leq \frac {L}{V}\\
\min\left\{ g(t-\frac{L}{V}) +H(\la(t)), C \right\}, & t>\frac{L}{V}}
\eqn
where the link queue size is 
\bqn
\la(t)&=& \cas{{ll} k_0 Vt-G(t), & t\leq \frac {L}{V} \\G(t-\frac{L}{V})-G(t)+k_0L, & t>\frac{L}{V}}
\eqn
\esq
Here the indicator function $H(z)$ for $z\geq 0$ is defined as
\bqs
H(z)&=&\lim_{\Delta t\to 0^+} \frac {z}{\Delta t}=\cas{{ll}0, & z=0\\+\infty, & z>0}
\eqs
Then we define the supply, $s(t)$, for the link between 0 and $L$: 
\bsq\label{def:supply}
\bqn
s(t)&=&\cas{{ll} \min\left\{(K-k_0) W +H(\gamma(t)), C\right\}, & t\leq \frac{L}{W}\\  \min\left\{ g(t-\frac{L}{W}) +H(\gamma(t)), C \right\}, &t>\frac{L}{W}}
\eqn
where the link vacancy size is 
\bqn
\gamma(t)&=& \cas{{ll} (K-k_0) W t-G(t), & t\leq \frac{L}{W}\\ G(t-\frac{L}{W})-G(t)+(K-k_0) L, &t>\frac{L}{W}}
\eqn
\esq
Finally, at the signalized intersection for $x=0$, we extend the following macroscopic junction model, which was first proposed in CTM \citep{daganzo1995ctm}, to calculate the boundary flow-rate:
\bqn
g(t)&=&\beta(t)  \min\{d(t),s(t)\}. \label{junctionmodel}
\eqn 
Thus from \refe{def:flow-rate}, \refe{def:demand}, \refe{def:supply}, and \refe{junctionmodel}, we obtain the continuous LTM for the signalized ring road, which is a delay-differential equation with $G(t)$ as the state variable. 

The discrete LTM with a time step-size of $\Delta t$, where $H(z)\dt=z$, can be calculated in the following steps:
\bsq \label{discrete-LTM}
(i) The discrete link demand is
\bqn
d(t)\dt &=&\cas{{ll} \min\{ (t+\dt)k_0 V-G(t), C \dt\}, & t+\dt\leq \frac{L}V\\ \min\{G(t+\dt-\frac LV)+k_0L-G(t), C\dt \}, & t+\dt>\frac LV}
\eqn
(ii) The discrete link supply is
\bqn
s(t) \dt &=&\cas{{ll}\min\{  (t+\dt)(K-k_0)W-G(t), C\dt \}, &t+\dt\leq \frac LW\\ \min\{G(t+\dt-\frac LW)+(K-k_0)L-G(t), C\dt\}, & t+\dt>\frac LW }
\eqn
(iii) The discrete boundary flow-rate is
\bqn
g(t)\dt &=& \beta(t) \min\{d(t)\dt, s(t) \dt\}.
\eqn
(iv) The boundary flow can be updated by
\bqn
G(t+\dt)&=&G(t)+g(t) \dt.
\eqn
\esq

From the discrete LTM, we can prove the following theorem, whose proof is given in Appendix A.
\begin{theorem} \label{thm:ltmsol} At a large time $t$, the continuous LTM for a signalized ring road can be solved by the following equation for the boundary flow: 
\bsq \label{continuousltm-large}
\bqn
G(t)&=&\min \{G(t-\frac LV)+k_0L,  G(t-\frac LW)+(K-k_0)L, G(iT)+(t-iT)C\}, \label{threewaves}
\eqn
for $t-iT\in (0, \pi T]$ during the effective green intervals, and 
\bqn
G(t)&=&G((i+\pi)T), \label{constantflow}
\eqn
 for $t-iT\in (\pi T, T]$ during the effective red intervals.
\esq
\end{theorem}

Note that, in \refe{threewaves}, three characteristic waves are considered when determining $G(t)$: the first one traveling forward at the free-flow speed, $V$; the second one traveling backward at the shock wave speed, $-W$; and the third one stationary at $x=0$. In contrast, when solving \refe{lwr_vt} inside the ring road with $x\in(0,L)$, only the first two characteristic waves need to be considered. This difference is caused by the existence of the traffic signal at $x=0$.

\section{Macroscopic fundamental diagram in stationary states}

In \citep{jin2014_hamiltonian}, it was shown that time-periodic solutions in both density and flow-rate exist for \refe{lwr_vt_signalized_ring} after a long time, and the period is the signal cycle length, $T$, with relatively large $T$, or multiple times of $T$ with relatively small $T$. In this section, we consider those with a period of $T$ as stationary states in the signalized ring road. We can see that the signalized ring road is in a stationary state if and only if $g(t+T)=g(t)$. In addition, stationary states can also be defined by  
\bqn
G(t+T)=G(t)+\bar g T, \label{stationarystate}
\eqn 
where $\bar g \in [0,\pi C]$ is the average flow-rate during a cycle.

In this section, we derive and analyze the macroscopic fundamental diagram (MFD), which is defined as the relationship between the average flow-rate, $\bar g$, and density, $k_0$, as well as signal settings in stationary states.

\subsection{Derivation of macroscopic fundamental diagram}
We divide the free-flow travel time, $\frac LV$, by the cycle length, $T$, to find the modulus, $j_1$, and the remainder, $\alpha_1$. That is,
\bsq \label{defmodrem}
\bqn
 \frac LV&=&\theta_1 T , \quad \theta_1=j_1+\alpha_1, \quad j_1=\lfloor \frac{L}{VT} \rfloor =0,1,\cdots, \quad 0\leq \alpha_1 <1,
\eqn 
where $\lfloor \cdot \rfloor$ is the floor function.
  Similarly we divide the shock wave propagation time, $\frac LW$, by the cycle length, $T$, to find the modulus, $j_2$, and the remainder, $\alpha_2$. That is
  \bqn
  \frac LW&=&\theta_2 T , \quad \theta_2=j_2+\alpha_2,  \quad j_2=\lfloor \frac{L}{WT} \rfloor =0,1,\cdots, \quad 0\leq \alpha_2 <1. 
\eqn
\esq

Further we define two critical densities, $k_1$ and $k_2$, as follows:
\bsq\label{twocriticaldensities}
\bqn
 k_1 &\equiv& \frac{j_1+\min\{\frac{\alpha_1}{\pi},1\}}{j_1+\alpha_1}  \pi \bar K,\\
 k_2 &\equiv& K-\frac{j_2+\min\{\frac{\alpha_2}{\pi},1\}}{j_2+\alpha_2} \pi \frac{C}{W}.
\eqn
\esq
Then we have the following Lemma.
\begin{lemma} \label{k1k2lemma1}$k_1$ and $k_2$ satisfy
\bqn
\pi \bar K \leq k_1 \leq \bar K \leq k_2 \leq K-\pi \frac {C}{W}. \label{k1k2relation}
\eqn
In addition, $k_1=\pi \bar K$ if and only if $\alpha_1=0$, and $k_1=\bar K$ if and only if $\pi T\geq \frac LV$; $k_2=K-\pi \frac CW$ if and only if $\alpha_2=0$, and $k_2=\bar K$ if and only if $\pi T\geq \frac LW$. 
\end{lemma}
{\em Proof}. From the definitions of $j_1$, $\alpha_1$, $j_2$, and $\alpha_2$ in \refe{defmodrem}, we can see that $\alpha_1\leq \frac{\alpha_1}\pi$ and $\alpha_1<1$. Thus $k_1\geq \pi \bar K$, where the equality holds if and only if $\alpha_1=0$. In addition, since $\pi<1$, $k_1\leq \frac{\pi j_1+\alpha_1}{j_1+\alpha_1} \bar K\leq \bar K$, where the equality holds if and only if $j_1=0$ and $\alpha_1\leq \pi$. That is, $k_1=\bar K$ if and only if $\pi T\geq \frac LV$. Similarly we can prove $\bar K\leq k_2 \leq K-\pi \frac CW$: $k_2=\bar K$ if and only if $\pi T\geq \frac LW$, and $k_2=K-\pi \frac CW$ if and only if $\alpha_2=0$. Thus \refe{k1k2relation} is correct. \eop

Then LTM during the effective green interval, \refe{threewaves}, leads to
\bqs
G(iT+\pi T)&=&\min \{G((i-j_1) T+(\pi-\alpha_1) T )+k_0L,  G((i-j_2)T +(\pi-\alpha_2) T)+(K-k_0)L, \nonumber \\&&G(iT)+ \pi T C\}=G(iT)+\bar g T. 
\eqs
LTM during the red interval, \refe{constantflow}, leads to $G((i+1)T)=G(iT+\pi T)$. Thus in stationary states, from \refe{stationarystate}  we have $G((i+1)T)=G(iT)+\bar g T$ and the following main equation for finding MFD:
\bqn
&&\min \{G((i-j_1) T+(\pi-\alpha_1) T )+k_0L,  G((i-j_2)T +(\pi-\alpha_2) T)+(K-k_0)L, \nonumber \\&&G(iT)+ \pi T C\}=G(iT)+\bar g T. \label{mainequation}
\eqn

Since $G(iT) +\bar g T\leq G(iT)+ \pi TC$, $\bar g \leq \pi C$. Clearly, the average flow-rate $\pi C$  can only be achieved when the boundary flow-rate is always maximum at capacity during the whole effective green interval; i.e., when $g(t)=C$ for $t-iT\in[0,\pi T]$.

Solving \refe{mainequation}, we obtain MFD for the signalized ring road in the following theorem, whose proof is given in Appendix B. 
\begin{theorem} \label{mfdtheorem} MFD for the signalized ring road is given by the following piecewise linear function:
\bqn
\bar g&=& \cas{{ll} \frac {k_0}{k_1} \pi C, & 0\leq k_0<k_1 \\ \pi C, & k_1\leq k_0 \leq k_2 \\ \frac {K-k_0}{K-k_2} \pi C, & k_2<k_0 \leq K } \label{mfd}
\eqn
\end{theorem}

\bfg\bc
\includegraphics[width=4in]{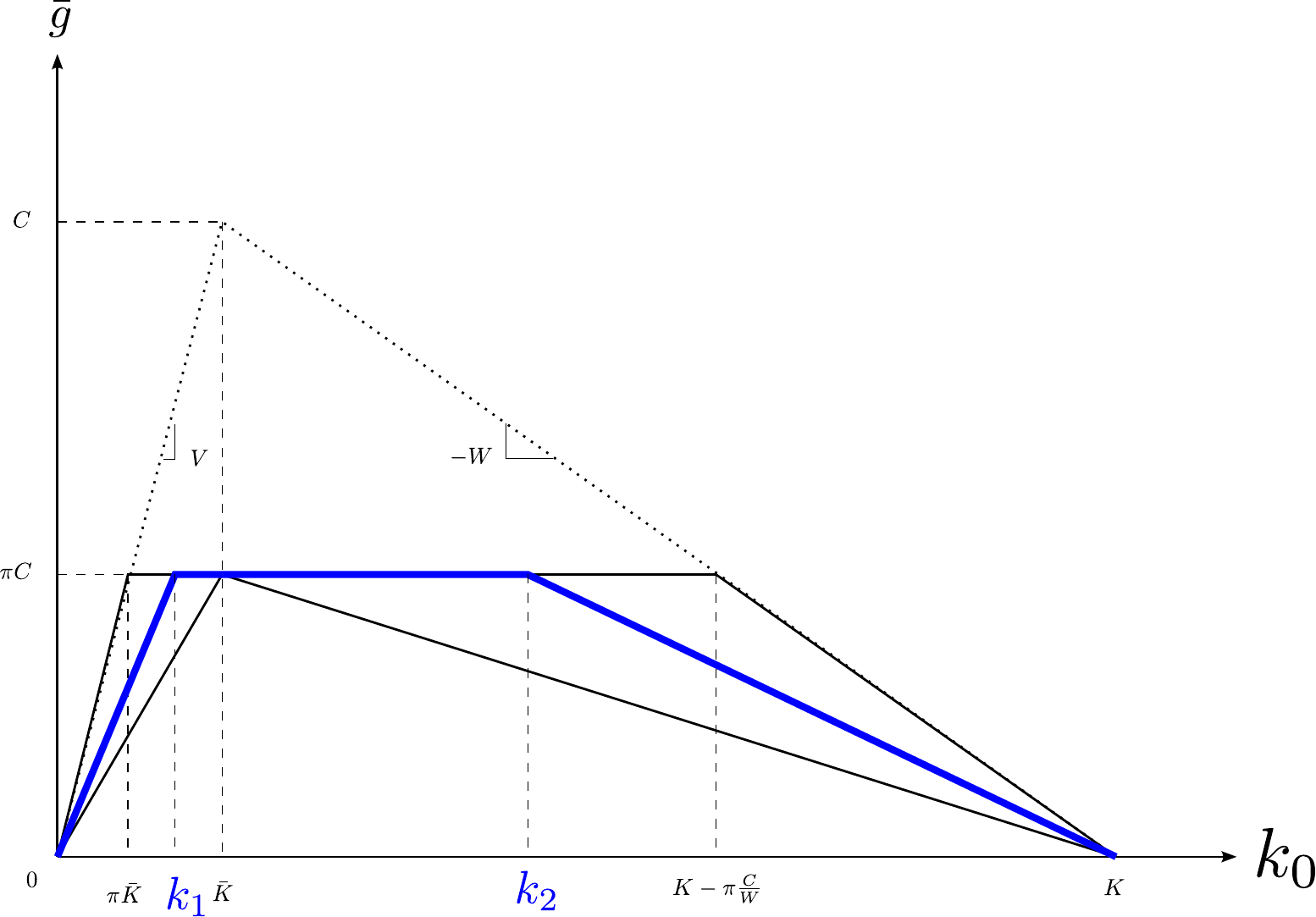}\caption{Macroscopic fundamental diagram for a signalized ring road}\label{tri_fd_signalized_ring}
\ec\efg

MFD in \refe{mfd} is shown in \reff{tri_fd_signalized_ring}, where the dotted lines are for the original triangular fundamental diagram, the thick solid lines for MFD, and the thin solid lines for the boundaries of MFD.
The derived MFD has the same shape as those in Figure 2 of \citep{gartner2004analysis}, which were obtained through simulations.
It is also consistent in principle with the piecewise linear MFD in \citep{daganzo2008analytical}, where a numerical method was provided to calculate MFD in a ring road with multiple signals. However, to the best of our knowledge, \refe{mfd} is the first explicit formula for MFD derived analytically. Such a formula is instrumental for further analysis of system performance and signal design.

\begin{corollary}
MFD for the signalized ring road can be re-written as 
\bqn
\bar g&=&\min\{\phi_1,\pi C, \phi_2\}, \label{mfd2}
\eqn
where $\phi_1=\frac{k_0}{k_1}\pi C$ and $\phi_2=\frac {K-k_0}{K-k_2} \pi C$.
\end{corollary}

\subsection{Properties of MFD}

Then from \refe{twocriticaldensities} we have the following properties of $\phi_1$ and $\phi_2$:
\ben
\item 
When $j_1=0$ and $0<\alpha_1<1$, 
\bsq
\bqn
\phi_1&=&\cas{{ll}\pi V k_0, & 0<\alpha_1 \leq \pi \\ \alpha_1 V k_0, & \pi<\alpha_1<1}
\eqn
When $j_1\geq 1$ and $0\leq \alpha_1<1$,
\bqn
\phi_1&=&\cas{{ll}\frac{j_1+\alpha_1}{\pi j_1+\alpha_1} \pi V k_0, & 0\leq \alpha_1 \leq \pi \\ \frac{j_1+\alpha_1}{j_1+1} Vk_0, & \pi<\alpha_1<1}
\eqn
\esq
Thus $\phi_1$ is continuous in $\theta_1$; $\phi_1$ retains the global minimum $\pi V k_0$ when $0<\theta_1=\alpha_1\leq \pi$,  reaches global maximum $V k_0$ when $\theta_1=j_1$, and reaches local minima $\frac{j_1+\pi}{j_1+1} Vk_0$ when $\theta_1=j_1+\pi$.

\item 
When $j_2=0$ and $0<\alpha_2<1$, 
\bsq
\bqn
\phi_2&=&\cas{{ll} \pi (K-k_0)  W, & 0<\alpha_2 \leq \pi \\  \alpha_2 (K-k_0) W, & \pi<\alpha_2<1}
\eqn
When $j_2\geq 1$ and $0\leq \alpha_2<1$,
\bqn
\phi_2&=&\cas{{ll}\frac{j_2+\alpha_2}{\pi j_2+\alpha_2} \pi (K-k_0) W, & 0\leq \alpha_2 \leq \pi \\ \frac{j_2+\alpha_2}{j_2+1} (K-k_0)W, & \pi<\alpha_2<1}
\eqn
\esq
Thus $\phi_2$ is continuous in $\theta_2$; $\phi_2$ retains the global minimum $ \pi (K-k_0) W$ when $0<\theta_2=\alpha_2\leq \pi$, reaches global maximum $(K-k_0) W$ when $\theta_2=j_2$, and reaches local minima $\frac{j_2+\pi}{j_2+1} (K-k_0) W$ when $\theta_2=j_2+\pi_2$.

\een

\commentout{
\begin{figure} \bc
$\ba{c}
\includegraphics[height=2in]{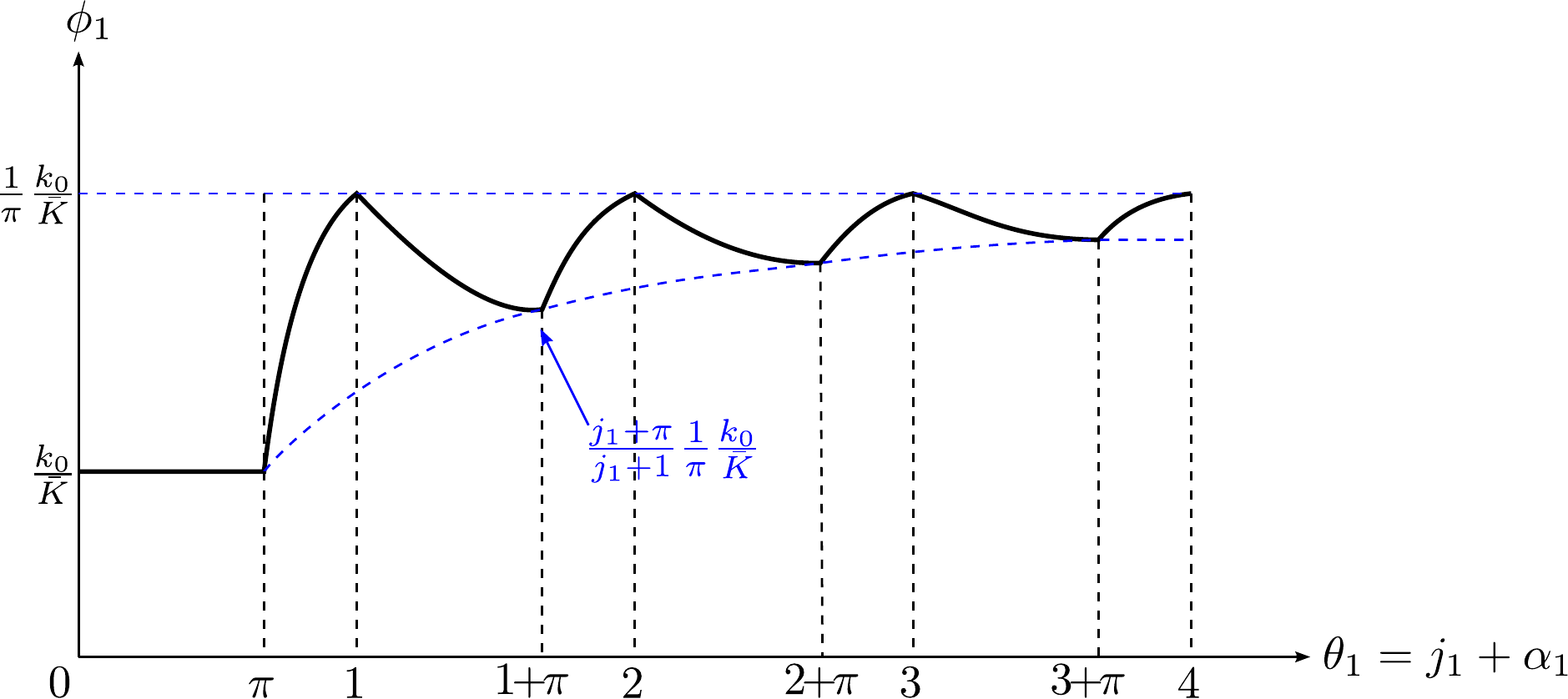} \\
\mbox{\bf (a)} \\
\includegraphics[height=2in]{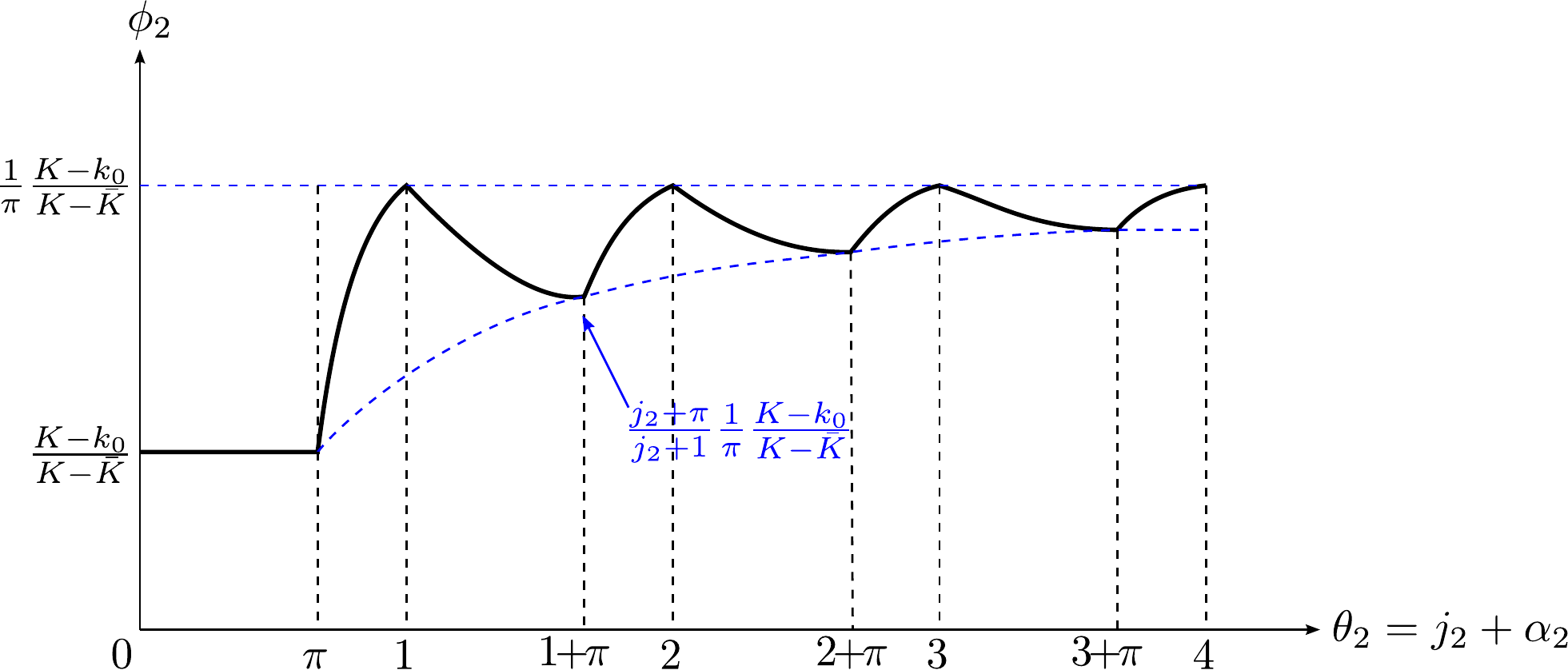} \\
    \mbox{\bf (b)}
\ea$
\caption{Illustrations of $\phi_1$ and $\phi_2$} \label{phi_theta} \ec 
\end{figure}
}

\begin{figure} \bc
$\ba{c}
\includegraphics[height=2in]{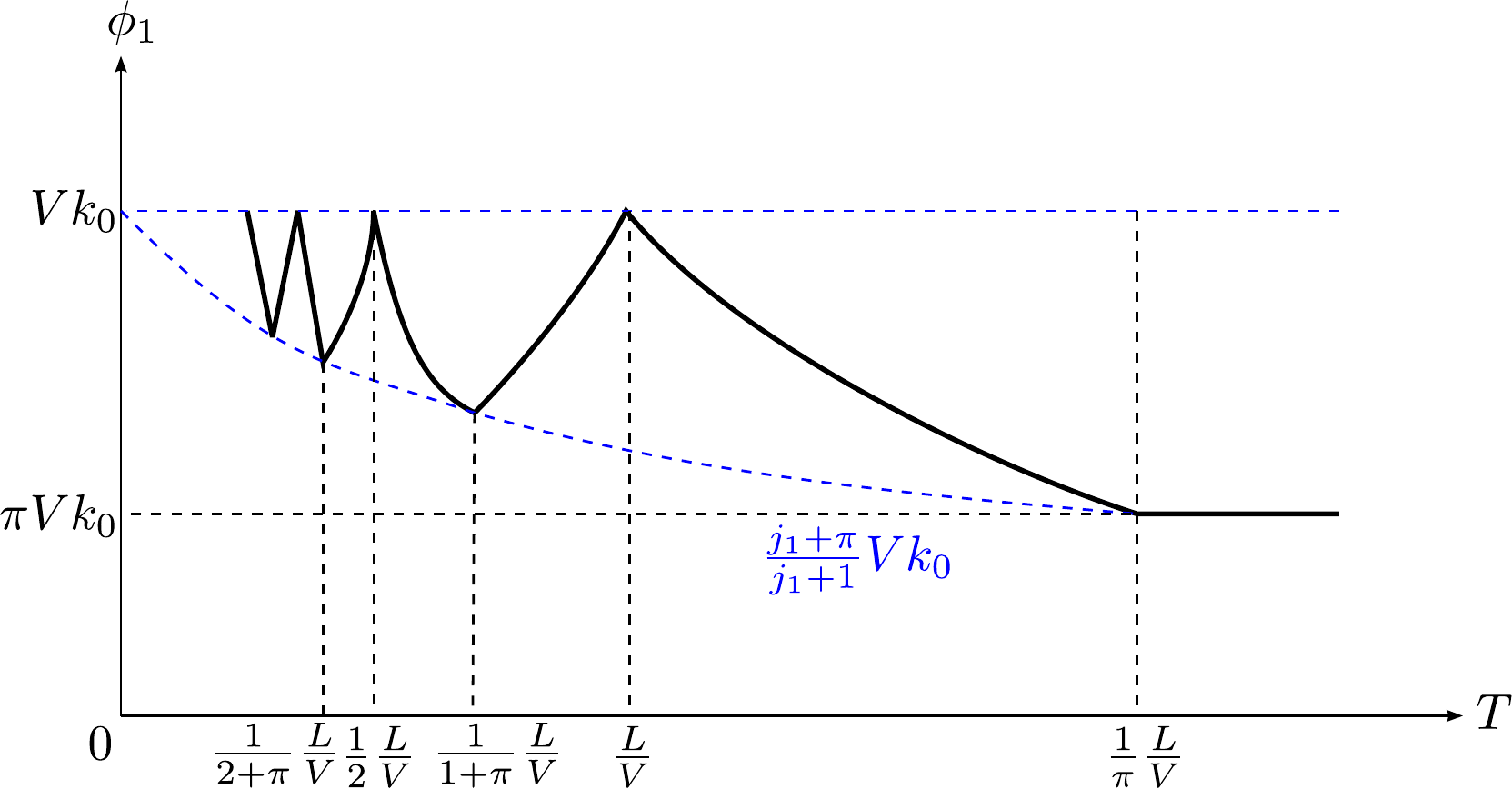} \\
\mbox{\bf (a)} \\
\includegraphics[height=2in]{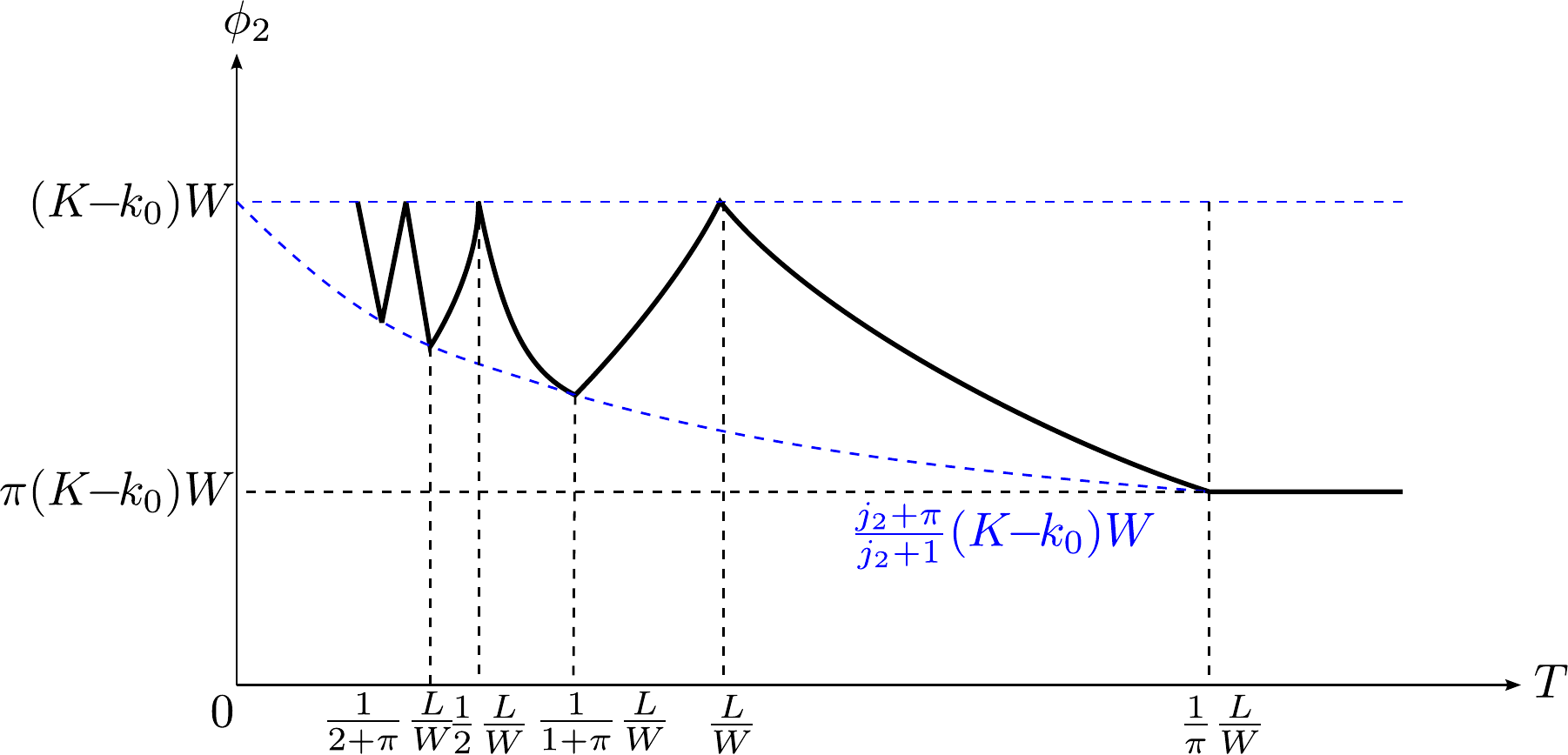} \\
    \mbox{\bf (b)}
\ea$
\caption{ (a) $\phi_1\sim T$ and (b) $\phi_2\sim T$ relations} \label{phi_T} \ec 
\end{figure}

Since $\theta_1=\frac LV\frac 1{T}$, and $\theta_2=\frac LW \frac 1{T}$, we can have the following $\phi_1\sim T$ and $\phi_2\sim T$ relations.
\begin{lemma}\label{lemmaphiT} $\phi_1$ and $\phi_2$ are functions of $T$, as shown in \reff{phi_T}:
\ben
\item $\phi_1$ is continuous in $T$; $\phi_1$ retains the global minimum $\pi V k_0$ when $T \geq \frac 1\pi \frac LV$, reaches global maximum $V k_0$ when $T=\frac 1{j_1} \frac LV$, and reaches local minima $\frac{j_1+\pi}{j_1+1} Vk_0$ when $T=\frac 1{j_1+\pi} \frac LV$.
 In particular, when $\frac LV \leq T \leq \frac 1\pi \frac LV$, the $\phi_1\sim T$ relation for the last decreasing branch is given by:
 \bqn
\phi_1= \frac{ k_0L}{T}. \label{phi1T}
\eqn

\item
$\phi_2$ is continuous in $T$; $\phi_2$ retains the global minimum $\pi (K-k_0) W$ when $T\geq \frac 1\pi \frac LW$, reaches global maximum $(K-k_0) W$  when $T=\frac 1{j_2} \frac LW$, and reaches local minima $\frac{j_2+\pi}{j_2+1} (K-k_0) W$ when $T=\frac 1{j_2+\pi} \frac LW$.
In particular, when $\frac LW \leq T \leq \frac 1\pi \frac LW$,  the $\phi_2\sim T$ relation for the last decreasing branch is given by:
\bqn
\phi_2= \frac{(K-k_0)L}{T}. \label{phi2T}
\eqn
\een
\end{lemma}

From \refe{mfd2}, we can see that the average flow-rate decreases in $k_1$ and increases in $k_2$. In particular, we have the following lemma.

\begin{lemma} \label{lemmaflowk1k2} In five regions for $k_0$, the average flow-rate varies with $k_1\in [\pi \bar K,\bar K]$ and $k_2\in[\bar K, K-\pi \frac CW]$ as follows:
\ben
\item When $k_0\in[0,\pi \bar K)$; i.e., traffic is very sparse, $\bar g=\phi_1$ for $k_1\in [\pi \bar K,\bar K]$, which decreases in $k_1$ and is independent of $k_2$. In this case, the global maximum flow-rate is $V k_0$ when $T=\frac 1{j_1} \frac LV$, and the global minimum flow-rate is $\pi V k_0$ when $T \geq \frac 1\pi \frac LV$.
\item When $k_0\in[\pi \bar K, \bar K)$; i.e., traffic is sparse, $\bar g= \min\{\phi_1 , \pi C\}$, which is first constant for $k_1\in [\pi \bar K, k_0]$ and then decreasing for $k_1 \in (k_0, \bar K]$, and independent of $k_2$. In this case, the global maximum flow-rate is $\pi C$, and the global minimum flow-rate is $\pi V k_0$ when $T \geq \frac 1\pi \frac LV$.
\item When $k_0=\bar K$; i.e., traffic is critical, $\bar g=\pi C$, which is constant for any $k_1$ and $k_2$.
\item When $k_0\in (\bar K, K-\pi \frac CW]$; i.e., traffic is dense, $\bar g=\min\{\phi_2,\pi C\}$, which is first increasing for $k_2\in[\bar K,k_0)$ and then constant for $k_2\in [k_0, K-\pi \frac CW]$, and independent of $k_1$. In this case, the global maximum flow-rate is $\pi C$, and the global minimum flow-rate is $\pi (K-k_0)W$ when $T \geq \frac 1\pi \frac LW$.

\item When $k_0\in (K-\pi \frac CW, K]$; i.e., traffic is very dense, $\bar g=\phi_2$ for $k_2\in[\bar K, K-\pi \frac CW]$, which increases in $k_2$ and is independent of $k_1$. In this case, the global maximum flow-rate is $(K-k_0)W$ when $T=\frac 1{j_2}\frac LW$, and the global minimum flow-rate is $\pi(K-k_0)W$ when $T \geq \frac 1\pi \frac LW$.
\een
\end{lemma}

From Lemmas \ref{lemmaphiT} and \ref{lemmaflowk1k2}, we can then determine $\bar g$ for any cycle length, $T$, at a given density, $k_0$. 
\commentout{
\begin{theorem}
At any density $k_0$, the average flow-rate, $\bar g$, is a continuous function of the cycle length, $T$, with the following properties: 
\ben
\item When $k_0\in[0,\pi \bar K)$, $\bar g$ retains the global minimum $\pi Vk_0$ when $T\geq \frac 1\pi \frac LV$, reaches global maximum $V k_0$ when $T=\frac 1{j_1} \frac LV$, and reaches local minima  $\frac{j_1+\pi}{j_1+1} V k_0$ when $T=\frac 1{j_1+\pi} \frac LV$. In particular, when $\frac LV<T<\frac 1\pi \frac LV$, $\bar g=\frac{k_0L}T$.
\item When $k_0\in[\pi \bar K, \bar K)$, $\bar g$ retains the global minimum $\pi Vk_0$ when $T\geq \frac 1\pi \frac LV$, reaches global maximum $V k_0$ when $T=\frac 1{j_1} \frac LV$, and reaches local minima  $\frac{j_1+\pi}{j_1+1} V k_0$ when $T=\frac 1{j_1+\pi} \frac LV$. But $\bar g$ is truncated by $\pi C$.
\item When $k_0=\bar K$, $\bar g=\pi C$, which is constant for any $T$.
\item When $k_0\in (\bar K, K-\pi \frac CW]$, $\bar g=\cas{{ll}\frac {K-k_0}{K-k_2} \pi C, & \bar K\leq k_2< k_0\\\pi C, &k_0\leq k_2\leq K-\pi \frac CW }$, which is first increasing and then constant in $k_2$, and independent of $k_1$. In this case, the maximum flow-rate is $\pi C$, and the minimum flow-rate is $\frac{K-k_0}{K-\bar K} \pi C$.

\item When $k_0\in (K-\pi \frac CW, K]$, $\bar g=\frac {K-k_0}{K-k_2} \pi C$ for $k_2\in[\bar K, K-\pi \frac CW]$, which increases in $k_2$ and is independent of $k_1$. In this case, the maximum flow-rate is $\pi C$, and the minimum flow-rate is $W(K-k_0)$.
\een
\end{theorem}
}

\section{Optimal signal settings}
In this section, we apply the analytical MFD formula, \refe{mfd}, to find optimal signal settings to  minimize the total travel time in stationary states, which is the time for all vehicles to traverse the whole road link and equals $\int_0^T A(0,t) -A(L,t) dt=\frac {(k_0 L)^2} {\bar g}$, since $A(0,t)=G(t)$, $A(L,t)=G(t)-k_0L$, the period equals $T=\frac {k_0 L}{\bar g}$. Here we obtain a delay formula for the stationary ring road,  $\frac {(k_0 L)^2} {\min\{\phi_1,\pi C,\phi_2\}}$, which is a function both density and signal settings.
At a given density $k_0$, therefore, the objective is equivalent to maximize the average flow-rate $\bar g$ and from \refe{mfd2} we have:
\bqn
\max \bar g=\max \min\{\phi_1,\pi C,\phi_2\}. \label{optimizationproblem1}
\eqn

If $\pi$ is constant and independent of $T$, we can see from Lemmas \ref{lemmaphiT} and \ref{lemmaflowk1k2} that there are an infinite number of solutions of $T$ for \refe{optimizationproblem1}: 
\ben
\item when $k_0< \bar K$, $\bar g$ reaches the global maximum, $\bar g^*=\min\{V k_0,\pi C\} $, when $T^*=\frac 1{j_1} \frac LV$ for any $j_1=1,2,\cdots$; 
\item when $k_0=\bar K$, any $T$ yield the same maximum flow-rate, $\bar g^*=\pi C$; 
\item when $k_0>\bar K$, $\bar g$ reaches the global maximum, $\bar g^*=\min\{(K- k_0)W,\pi C\} $,  when $T^*=\frac 1{j_2} \frac LW$ for any $j_2=1,2,\cdots$.
\een
Since $\bar g\geq \frac{j_1+\pi}{j_1+1} Vk_0$ for $k_0<\bar K$ and $\bar g\geq \frac{j_2+\pi}{j_2+1} (K-k_0)W$  for $k_0>\bar K$, $\bar g\to \bar g^*$ when $T\to 0$. That is, we can set the cycle length to be very small to achieve the best performance. In other words, it is best to install stop signs, which correspond to very small cycle lengths.

In reality, however, due to limited reaction times and bounded acceleration rates of drivers and vehicles, there exists a start-up lost time, and $\pi$ depends on $T$.

\subsection{Analytical solutions with a start-up lost time}
We denote the start-up lost time  in each phase by  $\delta$. Then the total effective green time for a cycle with two phases is only $T-2\delta$.
We assume that the effective green ratio is $\pi_0$, which allocates the total effective green time $T-2\delta$ to the ring road. Then the effective green time is $\pi T=(T-2\delta) \pi_0$. Therefore we have the following time-dependent ratio
\bqn
\pi &=& (1-\frac{2\delta}T) \pi_0. \label{defpi}
\eqn 
Therefore the average flow-rate $\bar g$ is bounded by
\bqn
\pi C&=&(1-\frac{2\delta}T) \pi_0 C, \label{upperbound}
\eqn
which equals 0 when $T=2\delta$ and increases in $T$.

Thus the objective function for the optimal control problem becomes
\bqn
\max \bar g=\max \min\{ \phi_1, (1-\frac{2\delta}T) \pi_0 C, \phi_2\},\label{optimizationproblem}
\eqn 
where $T\geq 2\delta$.
From Lemmas \ref{lemmaphiT} and \ref{lemmaflowk1k2}, we can see that in general $\phi_1$ and $\phi_2$ decreases in $T$ when $\pi$ is constant. Therefore it is possible to find a best cycle length to maximize the average flow-rate.

The average flow-rate, $\bar g= \min\{ \phi_1, (1-\frac{2\delta}T) \pi_0 C, \phi_2\}$, is a function of $L$, $K$, $V$, $W$, $\pi_0$, $\delta$, $T$, and $k_0$. Among these parameters, $L$ is determined by the road design and layout, $K$ by vehicle lengths and drivers' safety margins, $V$ by speed limit, $W$ by drivers' following aggressiveness, $\pi_0$ by effective green time allocation, $\delta$ by drivers' reaction times and vehicles' acceleration rates when a light turns green, $T$ by signal lengths, and $k_0$ by demand levels. In this study, we only consider the impacts of the three signal related parameters, $\pi_0$, $\delta$, and $T$.

First, \refe{mfd} can be written as
\bqn
\bar g&=& \cas{{ll} \frac {k_0}{j_1+\min\{\frac {\alpha_1}{\pi},1\}} \frac{L}{T}, & 0\leq k_0<k_1 \\
 \pi C, & k_1\leq k_0 \leq k_2 \\ 
 \frac {K-k_0}{j_2+\min\{\frac {\alpha_2}{\pi},1\}} \frac{L}{T}, & k_2<k_0 \leq K } 
\eqn
from which we can see that a larger $\pi$ leads to higher $\bar g$. Further from \refe{defpi} we can see that $\pi$ increases in $\pi_0$ and decreases in $\delta$. Thus, without changing other parameters, if we reduce the lost time $\delta$ and increase the effective green ratio $\pi_0$, we can increase the average flow-rate. The lost time, $\delta$, could be reduced by introducing autonomous or connected vehicles, which can have faster responses due to advanced sensors or communications between vehicles and traffic lights. However, the choice of $\pi_0$ has to be determined by the demand levels of other competing movements and is assumed to be constant. Thus here we assume that both $\delta$ and $\pi_0$ are fixed and aim to find an optimal cycle length.

When $\frac LV$ and $\frac LW$ are much larger than $2\delta$ (e.g., 10 times), $\pi\approx \pi_0$. Thus we assume that $\pi$ equals $\pi_0$ in $\phi_1$ and $\phi_2$. But note that the start-up lost time still impacts the average flow-rate in the MFD, as $\pi$ still depends on $T$ in \refe{upperbound}. Then solutions to the optimization problem, \refe{optimizationproblem}, are given by the following theorem.

\begin{theorem}\label{thm:optimalcycle} The optimal cycle lengths at different traffic densities are in the following:
\bsq
\ben
\item When traffic is very sparse with $k_0\in [0,\pi_0 \bar K)$, the maximum $\bar g^*\approx V k_0$, for which there exist multiple optimal cycle lengths:
\bqn
T^*= \frac 1{j_1} \frac LV,
\eqn
 for $j_1=1,2,\cdots$ and $Vk_0 \leq (1-\frac{2\delta}{T^*})\pi_0 C$.

\item When traffic is sparse with $k_0\in [\pi_0 \bar K,\bar K)$, the maximum $\bar g$ is determined by the intersection between $\pi C$ and the last decreasing branch of $\phi_1$ described by \refe{phi1T}:
\bqs
\bar g^*\approx\max_{T\in [\frac LV,\frac 1{\pi_0} \frac LV]} \min\{\frac {k_0 L}T, (1-\frac{2\delta}T) \pi_0 C\},
\eqs
for which there exists a unique optimal cycle length
\bqn
T^*&=&\frac{k_0 L}{\pi_0 C} +2\delta.
\eqn

\item When traffic is critical with $k_0=\bar K$, the maximum $\bar g$ is determined by $\pi C$:
\bqs
\bar g^*\approx\max_{T} \min\{(1-\frac{2\delta}{T})\pi_0 C \},
\eqs
fow which there exists a unique optimal cycle length 
\bqn
T^*&=&\infty.
\eqn

\item When traffic is dense with $k_0\in (\bar K, K-\pi \frac CW]$,  the maximum $\bar g$ is determined by the intersection between $\pi C$ and the last decreasing branch of $\phi_2$ described by \refe{phi2T}:
\bqs
\bar g^*\approx\max_{T\in [\frac LW,\frac 1{\pi_0} \frac LW]} \min\{\frac {(K-k_0) L}T, (1-\frac{2\delta}T) \pi_0 C\},
\eqs
for which there also exists a unique optimal cycle length
\bqn
T^*&=&\frac{(K-k_0) L}{\pi_0 C} +2\delta.
\eqn

\item When traffic is very dense with $k_0\in (K-\pi \frac CW,K]$, the maximum $\bar g^*\approx(K-k_0)W$, for which there exist multiple optimal cycle lengths:
\bqn
T^*=\frac 1{j_2} \frac LW,
\eqn
for $j_2=1,2,\cdots$ and $(K-k_0)W\leq (1-\frac{2\delta}{T^*})\pi_0 C$.

\een
\esq
\end{theorem}
{\em Proof}. The proof is straightforward, based on Lemmas \ref{lemmaphiT} and \ref{lemmaflowk1k2} as well as \refe{defpi}, and thus omitted. \eop

If we denote the congestion level by 
\bqn
\chi&=&\frac{\min\{Vk_0,C\}}{\min\{C,(K-k_0)W\}},
\eqn
which is the ratio of stationary demand over supply, then we have the following corollary from Theorem \ref{thm:optimalcycle}.

\begin{corollary} The optimal cycle lengths at different congestion levels are in the following:
\bsq\label{optimalcyclelength}
\ben
\item When traffic is very sparse with $\chi\in [0,\pi_0)$, the multiple optimal cycle lengths are
\bqs
T^*=\frac 1{j_1} \frac LV,
\eqs
 for $j_1=1,2,\cdots$ and $Vk_0 \leq (1-\frac{2\delta}{T^*})\pi_0 C$.

\item When traffic is sparse with $\chi\in [\pi_0,1)$, there exists a unique optimal cycle length
\bqn
T^*&=&\chi \frac{L}{\pi_0 V} +2\delta, \label{below1}
\eqn

\item When traffic is critical with $\chi=1$, there exists a unique optimal cycle length 
\bqs
T^*&=&\infty.
\eqs

\item When traffic is dense with $\chi \in (1, \frac 1{\pi_0}]$,   there exists a unique optimal cycle length
\bqn
T^*&=&\frac 1\chi \frac{L}{\pi_0 W} +2\delta. \label{above1}
\eqn

\item When traffic is very dense with $\chi \in (\frac 1{\pi_0},\infty)$, there exist multiple optimal cycle lengths:
\bqs
T^*=\frac 1{j_2} \frac LW,
\eqs
for $j_2=1,2,\cdots$ and $(K-k_0)W\leq (1-\frac{2\delta}{T^*})\pi_0 C$.

\een
\esq

\end{corollary}

When $\chi<1$, traditionally Webster's formula has been used to find the optimal cycle length \citep{roess2010traffic}: even though \refe{below1} is substantially different from Webster's optimal cycle length formula, it is consistent in principle with the latter, as it increases in both the congestion level and the lost time. But here we also obtain a simple formula \refe{above1} when $\chi>1$, and the optimal cycle length still increases in the lost time but decreases in the congestion level. In addition, the new formulas are derived from the LWR model and, therefore, more realistic.

\subsection{Numerical examples}
For a signalized ring road we choose the following parameters: $L=1200$ m, $V=20$ m/s, $W=5$ m/s, $K=1/7$ veh/m,  $\delta=3$ s, and $\pi_0=\frac 12$. Then $\pi_0 \bar K=\frac 12 \bar K$, $K-\pi_0 \frac CW=3 \bar K$, and $K=5\bar K$.
 
In \reff{optimal_cycle_length} we show the average flow-rates, calculated from \refe{mfd}, for different cycle lengths and densities. 
As shown in \reff{optimal_cycle_length}(a), when traffic is very sparse ($k_0<\pi_0 \bar K$), there exist multiple optimal cycle lengths, including $\frac LV$ and $\frac 12 \frac LV$. As shown in \reff{optimal_cycle_length}(b), when traffic is sparse ($k_0\in[\pi_0 \bar K, \bar K)$), there exists a unique optimal cycle length greater than $\frac LV$ but smaller than $2 \frac LV$: $T^*=\frac{k_0 L}{\pi_0 C}+2\delta$; in this case, $T=\frac LV$ still leads to near-optimal performance, but $T=2\frac LV$ is not acceptable, since the average flow-rate reaches the global minimum when $T=\frac 1{\pi_0} \frac LV+2\delta$. As shown in \reff{optimal_cycle_length}(c),  when traffic is dense ($k_0\in (\bar K, K-\pi_0 \frac CW]$), there exists a unique optimal cycle length greater than $\frac LW$ but smaller than $2\frac LW$: $T^*=\frac{(K-k_0)L}{\pi_0 C}+2\delta$; in this case, $T=\frac LW$ and even $T=\frac 12 \frac LW$ would lead to near-optimal performance. From \reff{optimal_cycle_length}(d), we can see that, when traffic is very dense ($k_0>K-\pi_0 \frac CW$), there can exist multiple optimal cycle lengths; in this case, $T=\frac 12 \frac LW$ would lead to near-optimal performance, and $T=\frac 14 \frac LW$ is also acceptable for very congested networks. These observations are consistent with the predictions in Theorem \ref{thm:optimalcycle}. In addition, from all the figures we can see that the maximum average flow-rate is unimodal in $k_0$; i.e., $\max \bar g$ increases in $k_0$ until $\bar K$ and then decreases.

\bfg\bc
\includegraphics[width=6.5in]{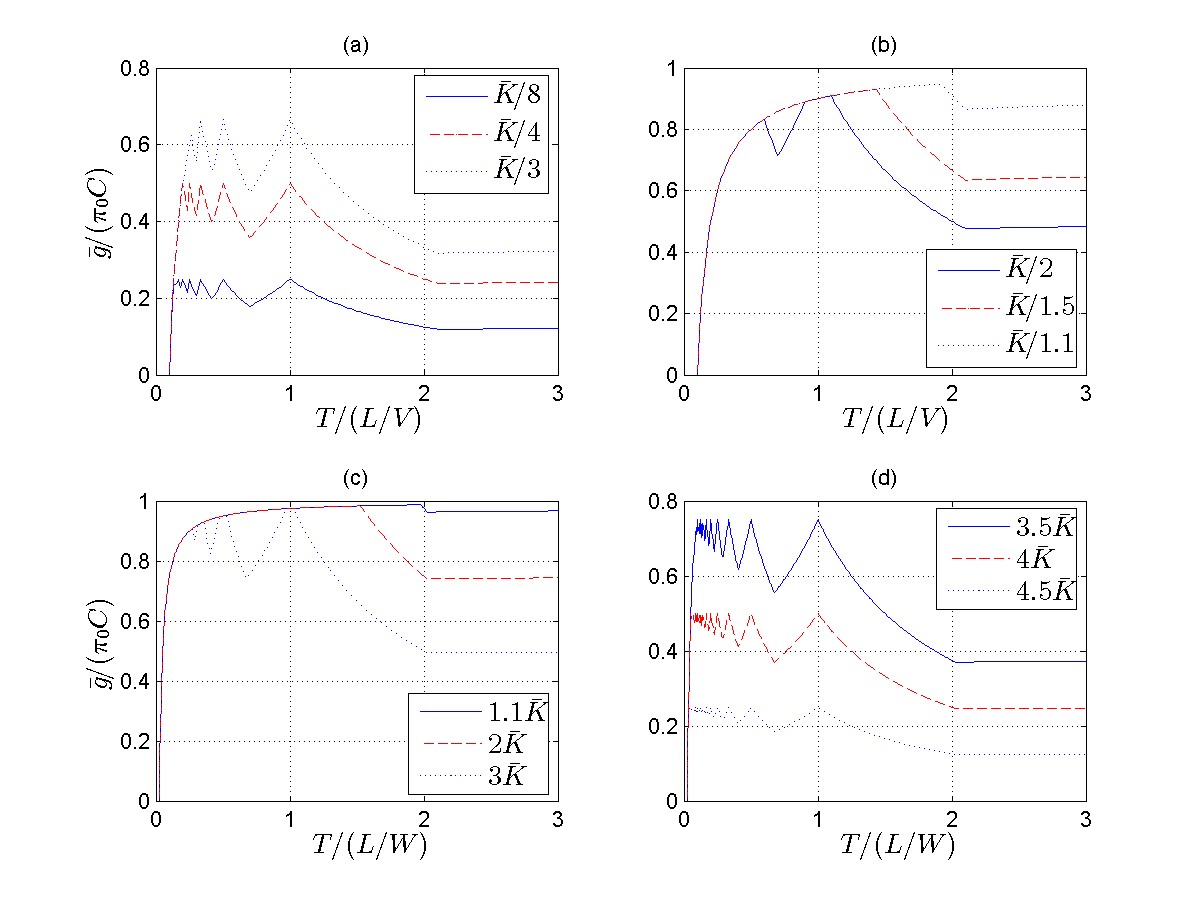} \caption{Average flow-rates vs cycle lengths for different densities with a start-up lost time: (a) very sparse traffic with $k_0 \in (0, \pi_0 \bar K)$; (b) sparse traffic with $k_0\in[\pi_0 \bar K, \bar K)$; (c) dense traffic with $k_0\in (\bar K, K-\pi_0 \frac CW]$; (d) very dense traffic with $k_0\in (K-\pi_0 \frac CW, K)$} \label{optimal_cycle_length}
\ec\efg

From \reff{optimal_cycle_length}(b) and \reff{optimal_cycle_length}(c), we can see that $\pi C$ slowly increases in $T$ when $T\geq \frac LV$. Therefore, we can choose a cycle length of $\frac LV=60$ s for sparse traffic and $\frac 12 \frac LW=120$ s for dense traffic. Such cycle lengths can lead to near-optimal performance. For examples, when $k_0=\bar K/1.5$, $\max \bar g=0.93 \pi_0 C$, and $\bar g=0.9 \pi_0 C$, which is about 3\% smaller, when $T=\frac LV$; when $k_0=2\bar K$, $\max \bar g=0.98 \pi_0 C$, and $\bar g=0.95 \pi_0 C$, which is also about 3\% smaller, when $T=\frac 12 \frac LW$. However, when $k_0=\bar K/1.5$, we cannot choose $T=2\frac LV=120$ s or larger, which leads to an average flow-rate of $0.67 \pi_0 C$, about 28\% smaller than the maximum value; similarly, when $k_0=2\bar K$, we cannot choose $T=2\frac LW=480$ s.
Furthermore, when traffic is extremely sparse with $k_0\leq \frac {\bar K}4$ shown in \reff{optimal_cycle_length}(a) or when traffic is very dense with $k_0>K-\pi_0 \frac CW$ shown in \reff{optimal_cycle_length}(d), stop signs, which can be considered as signals with very short cycle lengths, will be as effective as signals.

\section{Conclusion}
In this paper, we solved the link transmission model to obtain an equation for the boundary flow on a ring road with a pretimed signal. We then defined stationary states as periodic solutions with the cycle length as a period and derived the macroscopic fundamental diagram, from which we can calculate the average flow-rate from density and cycle length. After analyzing the impacts of the cycle length on the average flow-rate, we analytically derived optimal cycle lengths to maximize the average flow-rate subject to realistic start-up lost times due to drivers' reaction and acceleration behaviors. With numerical simulations, we verified the optimal solutions and suggested near-optimal cycle lengths under different congestion levels.

For the simplest signalized network, this study successfully fills the gap between methods based on delay formulas and those based on traffic simulation by presenting a new method that is both physically realistic and mathematically tractable. There are three particular contributions in this study. First, we obtained a simple link transmission model for the boundary flows on a signalized ring road, \refe{continuousltm-large}, which forms the foundation for solving and analyzing stationary states. Second, we derived an explicit macroscopic fundamental diagram, \refe{mfd}, in which the average flow-rate is a function of both traffic density and signal settings. Third, we presented formulas for optimal cycle lengths, \refe{optimalcyclelength}, under five levels of congestion with a start-up lost time. 

Even though it has been verified that the analytical optimal cycle length is quite accurate for large free-flow travel and shock wave propagation times, it is still important to check the results for really short links, where these times are small compared with the lost time. In the numerical solutions in Section 4.2, we found that some near-optimal cycle lengths can be used to avoid really long cycle lengths when the road is congested. But we will be interested in rigorously discussing such an optimization problem.

 This study for the simplest signalized network is a starting point for developing a unified framework for analyzing and designing signals in other networks, including streets with heterogeneous links and intersections, where cycle lengths, effective green ratios, offsets, and speed limits may be optimized simultaneously.
In addition, in this study, we analyze the performance of an open-loop control system and derive optimal cycle lengths for pretimed signals on a stationary signalized ring road, in which the plant is the whole traffic system, and the actuation is the signal. In the future, we will be interested in introducing feedback control mechanism to develop optimal signal settings under dynamic traffic conditions with random disturbances in general road networks.

\section*{Appendix A. Proof of Theorem \ref{thm:ltmsol}}
{\em Proof}. At a large time $t>\frac{L}{W}$, the discrete LTM, \refe{discrete-LTM}, can be written as 
\bqs
G(t+\dt)&=&G(t)+\beta(t) \min \{C\dt, \nonumber\\&&G(t+\dt-\frac LV)+k_0L-G(t), G(t+\dt-\frac LW)+(K-k_0)L-G(t) \} . 
\eqs
Then \refe{constantflow} is obvious for $t-iT\in (\pi T, T]$. In the following we prove that \refe{threewaves} is true for $t-iT\in (0, \pi T]$.

For $t-iT\in(0, \pi T]$ during the green intervals, $\beta(t)=1$, and
\bqs
G(t+\dt)&=&\min \{G(t+\dt-\frac LV)+k_0L, G(t+\dt-\frac LW)+(K-k_0)L, C\dt+G(t)\} . 
\eqs
 We denote $\dt=\frac{\pi T}{n}$. For $j=1$, we have 
\bqs
G(iT+\dt)&=&\min \{G(iT+\dt-\frac LV)+k_0L, G(iT+\dt-\frac LW)+(K-k_0)L, C\dt+G(iT)\}.
\eqs
We assume for any $1\leq j<n$,
\bqs
G(iT+j\dt)&=&\min \{G(iT+j\dt-\frac LV)+k_0L, G(iT+j\dt-\frac LW)+(K-k_0)L, Cj\dt+G(iT)\}.
\eqs 
Then for $j+1$, we have
\bqs
G(iT+(j+1)\dt)&=&\min \{G(iT+(j+1)\dt-\frac LV)+k_0L, G(iT+(j+1)\dt-\frac LW)+(K-k_0)L,\\
&& C\dt+G(iT+j\dt)\},\\
&=&\min \{G(iT+(j+1)\dt-\frac LV)+k_0L, G(iT+(j+1)\dt-\frac LW)+(K-k_0)L,\\
&& C\dt+G(iT+j\dt-\frac LV)+k_0L, C\dt+G(iT+j\dt-\frac LW)+(K-k_0)L,\\
&& C\dt+Cj\dt+G(iT)\}.
\eqs
Since $G(iT+(j+1)\dt-\frac LV)\leq C\dt+G(iT+j\dt-\frac LV)$ and $G(iT+(j+1)\dt-\frac LW)\leq C\dt+G(iT+j\dt-\frac LW)$, we have
\bqs
G(iT+(j+1)\dt)&=&\min \{G(iT+(j+1)\dt-\frac LV)+k_0L, G(iT+(j+1)\dt-\frac LW)+(K-k_0)L,\\
&& C(j+1)\dt+G(iT)\}.
\eqs
Thus from the method of induction, we have for $j=1,\cdots,n$
\bqs
G(iT+j\dt)&=&\min \{G(iT+j\dt-\frac LV)+k_0L, G(iT+j\dt-\frac LW)+(K-k_0)L, Cj\dt+G(iT)\}.
\eqs
If we denote $t=iT+j \dt$ ($1\leq j\leq n$), then we obtain \refe{threewaves}. 
\eop

\section*{Appendix B. Proof of Theorem \ref{mfdtheorem}} 
{\em Proof}. We derive \refe{mfd} in the following three cases.
\ben
\item 
From \refe{mainequation}, we can see that $\bar g=\pi C$ if and only if
\bqs
G((i-j_1) T+(\pi-\alpha_1) T )+k_0L &\geq & G(iT)+ \pi T C ,\\
G((i-j_2)T +(\pi-\alpha_2) T)+(K-k_0)L & \geq &G(iT)+ \pi T C.
\eqs
For the first equation, we have the following two scenarios: 
\ben
\item If $\pi \leq \alpha_1$, then $i-j_1 -1+\pi <i-j_1+\pi -\alpha_1\leq i-j_1$, and $G(iT)-G((i-j_1) T+(\pi-\alpha_1) T )=j_1 \bar g T =j_1 \pi C T$. Thus $k_0 \geq (j_1+1) \frac{\pi TC} L=\frac{j_1+1}{j_1+\alpha_1} \pi \bar K$.
\item If $\pi >\alpha_1$, then $i-j_1<i-j_1+\pi-\alpha_1\leq i-j_1+\pi$, and $G(iT)-G((i-j_1) T+(\pi-\alpha_1) T )=j_1 \pi C T -(\pi-\alpha_1)C T$. Thus $k_0 \geq (j_1+\frac{\alpha_1}{\pi}) \frac{\pi TC} L =\frac{j_1+\frac{\alpha_1}{\pi}}{j_1+\alpha_1} \pi \bar K$.
\een
Thus the first equation is equivalent to $k_0\geq k_1$. Similarly we can prove that the second equation is equivalent to $k_0\leq k_2$. Therefore the second scenario in \refe{mfd} is proved.

\item
When $G((i-j_1) T+(\pi-\alpha_1) T )+k_0L < G(iT)+ \pi T C$; i.e., if $k_0<k_1\leq k_2$, then $\bar g<\pi C$, and from \refe{mainequation} we have
\bqs
G(iT)+\bar g T&=& G((i-j_1) T+(\pi-\alpha_1) T )+k_0L.
\eqs
When $\pi\leq \alpha_1$, we have $\bar g=\frac{k_0 L}{(j_1+1)T}=\frac {k_0}{k_1} \pi C$. When $\pi>\alpha_1$, we have
\bqs
(j_1+1)T \bar g-\int_0^{(\pi-\alpha_1)T} g(t) dt&=&k_0 L. 
\eqs
If we assume that $g(t)$ is evenly distributed between $0$ an $\pi T$, then  $\int_0^{(\pi-\alpha_1)T} g(t) dt=(1-\frac{\alpha_1}\pi) T \bar g$, and $\bar g=\frac {k_0}{k_1} \pi C$. Therefore the first scenario in \refe{mfd} is proved.

\item
 When $G((i-j_2)T +(\pi-\alpha_2) T)+(K-k_0)L< G(iT)+ \pi T C$; i.e., if $k_0>k_2$, then $\bar g<\pi C$, and
from \refe{mainequation} we have
\bqs
G(iT)+\bar g T&=& G((i-j_2)T +(\pi-\alpha_2) T)+(K-k_0)L.
\eqs
When $\pi\leq \alpha_2$, we have $\bar g=\frac{(K-k_0) L}{(j_2+1)T}=\frac {K-k_0}{K-k_2} \pi C$. When $\pi>\alpha_2$, we have
\bqs
(j_2+1)T \bar g-\int_0^{(\pi-\alpha_2)T} g(t) dt&=&(K-k_0) L. 
\eqs
If we assume that $g(t)$ evenly distributes from $0$ to $\pi T$, then  $\int_0^{(\pi-\alpha_2)T} g(t) dt=(1-\frac{\alpha_2}\pi) T \bar g$, and $\bar g=\frac {K-k_0}{K-k_2} \pi C$. Therefore the third scenario in \refe{mfd} is proved.
\een
\eop

\section*{Acknowledgments}
The first author would like to thank the support of a UCCONNECT faculty grant.

\section*{References}
\pdfbookmark[1]{References}{references}


\begin{thebibliography}{42}
\expandafter\ifx\csname natexlab\endcsname\relax\def\natexlab#1{#1}\fi
\expandafter\ifx\csname url\endcsname\relax
  \def\url#1{\texttt{#1}}\fi
\expandafter\ifx\csname urlprefix\endcsname\relax\fi

\bibitem[{Ardekani and Herman(1987)}]{ardekani1987urban}
Ardekani, S., Herman, R., 1987. Urban network-wide traffic variables and their
  relations. Transportation Science 21~(1), 1--16.

\bibitem[{Buisson and Ladier(2009)}]{buisson2009exploring}
Buisson, C., Ladier, C., 2009. Exploring the impact of homogeneity of traffic
  measurements on the existence of macroscopic fundamental diagrams.
  Transportation Research Record: Journal of the Transportation Research Board
  2124, 127--136.

\bibitem[{Cassidy et~al.(2011)Cassidy, Jang, and
  Daganzo}]{cassidy2011macroscopic}
Cassidy, M., Jang, K., Daganzo, C.~F., 2011. Macroscopic fundamental diagrams
  for freeway networks. Transportation Research Record: Journal of the
  Transportation Research Board 2260, 8--15.

\bibitem[{Chang and Lin(2000)}]{chang2000optimal}
Chang, T.-H., Lin, J.-T., 2000. Optimal signal timing for an oversaturated
  intersection. Transportation Research Part B: Methodological 34~(6),
  471--491.

\bibitem[{Chang and Sun(2004)}]{chang2004modeling}
Chang, T.-H., Sun, G.-Y., 2004. Modeling and optimization of an oversaturated
  signalized network. Transportation Research Part B 38~(8), 687--707.

\bibitem[{Daganzo(1995)}]{daganzo1995ctm}
Daganzo, C.~F., 1995. {The cell transmission model {II}: Network traffic}.
  Transportation Research Part B 29~(2), 79--93.

\bibitem[{Daganzo(2005{\natexlab{a}})}]{daganzo2005variationalKW}
Daganzo, C.~F., 2005{\natexlab{a}}. A variational formulation of kinematic
  waves: basic theory and complex boundary conditions. Transportation Research
  Part B 39~(2), 187--196.

\bibitem[{Daganzo(2005{\natexlab{b}})}]{daganzo2005variationalKW2}
Daganzo, C.~F., 2005{\natexlab{b}}. A variational formulation of kinematic
  waves: Solution methods. Transportation Research Part B 39~(10), 934--950.

\bibitem[{Daganzo(2007)}]{daganzo2007gridlock}
Daganzo, C.~F., 2007. {Urban gridlock: Macroscopic modeling and mitigation
  approaches}. Transportation Research Part B 41~(1), 49--62.

\bibitem[{Daganzo et~al.(2011)Daganzo, Gayah, and
  Gonzales}]{daganzo2011bifurcations}
Daganzo, C.~F., Gayah, V.~V., Gonzales, E.~J., 2011. Macroscopic relations of
  urban traffic variables: Bifurcations, multivaluedness and instability.
  Transportation Research Part B 45~(1), 278--288.

\bibitem[{Daganzo and Geroliminis(2008)}]{daganzo2008analytical}
Daganzo, C.~F., Geroliminis, N., 2008. An analytical approximation for the
  macroscopic fundamental diagram of urban traffic. Transportation Research
  Part B 42~(9), 771--781.

\bibitem[{D'ans and Gazis(1976)}]{dans1976optimal}
D'ans, G., Gazis, D., 1976. Optimal control of oversaturated store-and-forward
  transportation networks. Transportation Science 10~(1), 1--19.

\bibitem[{Dion et~al.(2004)Dion, Rakha, and Kang}]{dion2004intersection}
Dion, F., Rakha, H., Kang, Y., 2004. {Comparison of delay estimates at
  under-saturated and over-saturated pre-timed signalized intersections}.
  Transportation Research Part B 38~(2), 99--122.

\bibitem[{Evans(1998)}]{evans1998pde}
Evans, L., 1998. {Partial Differential Equations}. American Mathematical
  Society.

\bibitem[{Gartner et~al.(1975)Gartner, Little, and
  Gabbay}]{gartner1975optimization}
Gartner, N.~H., Little, J.~D., Gabbay, H., 1975. Optimization of traffic signal
  settings by mixed-integer linear programming: Part i: The network
  coordination problem. Transportation Science 9~(4), 321--343.

\bibitem[{Gartner and Wagner(2004)}]{gartner2004analysis}
Gartner, N.~H., Wagner, P., 2004. Analysis of traffic flow characteristics on
  signalized arterials. Transportation Research Record: Journal of the
  Transportation Research Board 1883, 94--100.

\bibitem[{Gazis and Potts(1963)}]{gazis1963over}
Gazis, D., Potts, R., 1963. The over-saturated intersection. In: Procedings of
  the 2nd International Symposium on Traffic Theory. Organisation for Economic
  Co-operation and Development, pp. 221--237.

\bibitem[{Gazis(1964)}]{gazis1964optimum}
Gazis, D.~C., 1964. Optimum control of a system of oversaturated intersections.
  Operations Research 12~(6), 815--831.

\bibitem[{Geroliminis and Boyaci(2012)}]{geroliminis2012effect}
Geroliminis, N., Boyaci, B., 2012. The effect of variability of urban systems
  characteristics in the network capacity. Transportation Research Part B
  46~(10), 1607--1623.

\bibitem[{Geroliminis and Daganzo(2008)}]{geroliminis2008eus}
Geroliminis, N., Daganzo, C.~F., 2008. {Existence of urban-scale macroscopic
  fundamental diagrams: Some experimental findings}. Transportation Research
  Part B 42~(9), 759--770.

\bibitem[{Geroliminis et~al.(2013)Geroliminis, Haddad, and
  Ramezani}]{geroliminis2013optimal}
Geroliminis, N., Haddad, J., Ramezani, M., 2013. Optimal perimeter control for
  two urban regions with macroscopic fundamental diagrams: A model predictive
  approach. Intelligent Transportation Systems, IEEE Transactions on 14~(1),
  348--359.

\bibitem[{Godfrey(1969)}]{godfrey1969mechanism}
Godfrey, J., 1969. The mechanism of a road network. Traffic Engineering and
  Control 8~(8), 323--327.

\bibitem[{Haberman(1977)}]{haberman1977model}
Haberman, R., 1977. Mathematical models. Prentice Hall, Englewood Cliffs, NJ.

\bibitem[{Improta and Cantarella(1984)}]{improta1984control}
Improta, G., Cantarella, G., 1984. Control system design for an individual
  signalized junction. Transportation Research Part B: Methodological 18~(2),
  147--167.

\bibitem[{Jin(2015)}]{jin2015ltm}
Jin, W.-L., 2015. Continuous formulations and analytical properties of the link
  transmission model. Transportation Research Part B 74, 88--103.

\bibitem[{Jin et~al.(2013)Jin, Gan, and Gayah}]{jin2013dring}
Jin, W.-L., Gan, Q.-J., Gayah, V.~V., 2013. {A kinematic wave approach to
  traffic statics and dynamics in a double-ring network}. Transportation
  Research Part B 57, 114--131.

\bibitem[{Jin and Yu(2015)}]{jin2014_hamiltonian}
Jin, W.-L., Yu, Y., 2015. {Asymptotic solution and effective Hamiltonian of a
  Hamilton-Jacobi equation in the modeling of traffic flow on a homogeneous
  signalized road}. Journal de Math\'ematiques Pures et Appliqu\'ees (JMPA)In
  press.

\bibitem[{Lighthill and Whitham(1955)}]{lighthill1955lwr}
Lighthill, M.~J., Whitham, G.~B., 1955. {On kinematic waves: II. A theory of
  traffic flow on long crowded roads}. Proceedings of the Royal Society of
  London A 229~(1178), 317--345.

\bibitem[{Mahmassani et~al.(1987)Mahmassani, Williams, and
  Herman}]{mahmassani1987performance}
Mahmassani, H., Williams, J., Herman, R., 1987. Performance of urban traffic
  networks. Proceedings of the Tenth International Symposium on Transportation
  and Traffic Theory.

\bibitem[{Miller(1963)}]{miller1963settings}
Miller, A.~J., 1963. Settings for fixed-cycle traffic signals. Operations
  Research 14~(4), 373--386.

\bibitem[{Moskowitz(1965)}]{moskowitz1965discussion}
Moskowitz, K., 1965. { Discussion of `freeway level of service as in uenced by
  volume and capacity characteristics' by D.R. Drew and C. J. Keese}. Highway
  Research Record 99, 43--44.

\bibitem[{Munjal et~al.(1971)Munjal, Hsu, and Lawrence}]{munjal1971multilane}
Munjal, P.~K., Hsu, Y.~S., Lawrence, R.~L., 1971. Analysis and validation of
  lane-drop effects of multilane freeways. Transportation Research 5~(4),
  257--266.

\bibitem[{Newell(1989)}]{newell1989theory}
Newell, G.~F., 1989. Theory of highway traffic signals. Tech. rep.

\bibitem[{Newell(1993)}]{newell1993sim}
Newell, G.~F., 1993. {A simplified theory of kinematic waves in highway traffic
  {I}: General theory. {II}: Queuing at freeway bottlenecks. {III}:
  Multi-destination flows}. Transportation Research Part B 27~(4), 281--313.

\bibitem[{Olszewski et~al.(1995)Olszewski, Fan, and Tan}]{olszewski1995area}
Olszewski, P., Fan, H., Tan, Y., 1995. Area-wide traffic speed-flow model for
  the singapore cbd. Transportation Research Part A 29~(4), 273--281.

\bibitem[{Papageorgiou(1995)}]{papageorgiou1995integrated}
Papageorgiou, M., 1995. An integrated control approach for traffic corridors.
  Transportation Research Part C: Emerging Technologies 3~(1), 19--30.

\bibitem[{Papageorgiou et~al.(2005)Papageorgiou, Diakaki, Dinopoulou,
  Kotsialos, and Wang}]{papageorgiou2005review}
Papageorgiou, M., Diakaki, C., Dinopoulou, V., Kotsialos, A., Wang, Y., 2005.
  {Review of road traffic control strategies}. Proceedings of the IEEE 91~(12),
  2043--2067.

\bibitem[{Park et~al.(1999)Park, Messer, and Urbanik}]{park1999traffic}
Park, B., Messer, C.~J., Urbanik, T., 1999. Traffic signal optimization program
  for oversaturated conditions: genetic algorithm approach. Transportation
  Research Record: Journal of the Transportation Research Board 1683, 133--142.

\bibitem[{Richards(1956)}]{richards1956lwr}
Richards, P.~I., 1956. Shock waves on the highway. Operations Research 4~(1),
  42--51.

\bibitem[{Roess et~al.(2010)Roess, Prassas, and McShane}]{roess2010traffic}
Roess, R., Prassas, E., McShane, W., 2010. {Traffic engineering}. Prentice
  Hall.

\bibitem[{Yperman(2007)}]{yperman2007link}
Yperman, I., 2007. The link transmission model for dynamic network loading.
  Ph.D. thesis.

\bibitem[{Yperman et~al.(2006)Yperman, Logghe, Tampere, and
  Immers}]{yperman2006mcl}
Yperman, I., Logghe, S., Tampere, C., Immers, B., 2006. {The Multi-Commodity
  Link Transmission Model for Dynamic Network Loading}. Proceedings of the TRB
  Annual Meeting.

\end{thebibliography}
\end{document}